\documentclass[12pt, 14paper,reqno]{amsart}
\usepackage{amsmath,amsfonts,amssymb}
\usepackage{mathrsfs}
\usepackage{longtable}
\usepackage[breaklinks]{hyperref}
\usepackage{graphicx}
\makeatletter
\@namedef{subjclassname@2020}{%
  \textup{2020} Mathematics Subject Classification}
\makeatother

\theoremstyle{plain}
\newtheorem{thm}{Theorem}[section]
\newtheorem{lem}{Lemma}[section]
\newtheorem{cor}{Corollary}[section]

\newtheorem{thma}{Theorem}

\theoremstyle{proof}

\numberwithin{equation}{section}

\begin{document} 
\title[Generalized Mersenne Numbers of the form $cx^2$]{Generalized Mersenne Numbers of the form $cx^2$}
\author{Azizul Hoque}
\email{ahoque.ms@gmail.com}

\subjclass[2020]{11D61; 11N32}

\date{\today}

\keywords{Generalized Mersene number, Diophantine equation, Integer solution}

\begin{abstract} Generalized Mersenne numbers are defined as $M_{p,n} = p^n- p + 1$, where $p$ is any prime and $n$ is any positive integer. Here, we prove that for each pair $(c, p)$ with $c\geq 1$ an integer, there is at most one $M_{p, n}$ of the form $cx^2$ with a few exceptions. 

 \end{abstract}
\maketitle{}

\section{Introduction}
Mersenne numbers are positive integers of the form $2^n-1$ with $n \geq  1$. These numbers have attracted a great deal of interest since the seventeenth century. Furthermore, the primes of this form; so called Mersenne primes are traceable back to Euclid, who in his `\textit{Elements}' connected primes of the form $2^n - 1$ to even perfect numbers. In particular, Euclid-Euler theorem states that an even number is perfect if and only if it has the form $2^{n-1}(2^n -1)$, where $2^n -1$ is a prime number.
A perfect number is a positive integer that is equal to the sum of its proper divisors. Several earliest results spawned from attempts to understand these numbers. Although some modern researchers continue to attribute the same mystical significance to these numbers that the ancient people once did, these numbers remain a substantial inspiration for research in number theory (see \cite{D14, FLS09, P86, W83}). One of the challenging unsolved problems in number theory is Lenstra-Pomerance-Wagstaff conjecture, which states that there are infinitely many Mersenne primes. 

In \cite{HS14}, the author and Saikia studied a generalization of Mersenne numbers; so called generalized Mersenne numbers. These numbers are defined as $M_{p, n}=p^n-p+1$, where $p$ is any prime and $n$ is any positive integer. In this case too, we expect an extended version of Lenstra-Pomerance-Wagstaff conjecture. The precise problem is \textit{whether there are infinitely many primes of the form $M_{p, n}$ for each prime $p$}. A weaker version of this problem was posted in \cite{HS14}. Here, we investigate the problem:
   \textit{How many generalized Mersenne numbers are there of the form $cx^2$ for each pair of integers $(c, p)$ with $c\geq 1$ an integer and $p$ a prime?} 
Precisely, we prove:  

\begin{thm}\label{thm}
For any odd integer $c\geq 1$ and a prime $p$, the generalized Mersenne numbers of the form $cx^2$ are $1, 25$ and $121$, with at most one more possibility for each pair $(c, p)$. Further for even integer $c\geq 2$, there is no generalized Mersenne number of the form $cx^2$. 
\end{thm}

Assume that $2^n-1=cx^2$. Then for $n\geq 3$, we have $c\equiv 7\pmod 8$ and thus by \cite[p. 1]{CO00}, $cx^2+1=2^n$ has at least one solution $(x, n)$. Therefore we have the following straightforward corollary.
\begin{cor}
Let $c$ be a positive integer. If $c\not\equiv 7\pmod 8$, then $1$ is the only Mersenne number of the form $cx^2$. Further for $c\equiv 7\pmod 8$, there is exactly one Mersenne number of the form $cx^2$. 
\end{cor}   

The proof of Theorem \ref{thm} largely relies on a remarkable result of Bugeaud and Shorey \cite[Theorem 1]{BS01} on the positive integer solutions of certain Diophantine equations. 

\section{Preliminary descent}
We begin this section with a classical result of Bugeaud and Shorey \cite{BS01} on the number of positive integer solutions of certain Diophantine equations. Before stating this result, we need to introduce some definitions and notations.

Let $F_k$ (resp. $L_k$) denote the $k$-th term in the Fibonacci (resp. Lucas) sequence defined by $F_0=0,   F_1= 1$,
and $F_{k+2}=F_k+F_{k+1}$ (resp. $L_0=2,  L_1=1$, and $L_{k+2}=L_k+L_{k+1}$), where $k\geq 0$ is an integer. 
Given $\lambda\in \{1, \sqrt{2}, 2\}$, we define the sets $\mathcal{F}, \ \mathcal{G},\ \mathcal{H}\subset \mathbb{N}\times\mathbb{N}\times\mathbb{N}$ as follows:
\begin{align*}
\mathcal{F}&:=\{(F_{k-2\varepsilon},L_{k+\varepsilon},F_k)\,|\,
k\geq 2,\varepsilon\in\{\pm 1\}\},\\
\mathcal{G}&:=\{(1,4p^r-1,p)\,|\,\text{$p$ is an odd prime},r\geq 1\},\\
\mathcal{H}&:=\left\{(D_1,D_2,p)\,\left|\,
\begin{aligned}
&\text{$D_1$, $D_2$ and $p$ are mutually coprime positive integers with $p$}\\
&\text{an odd prime and there exist positive integers $r$, $s$ such that}\\
&\text{$D_1s^2+D_2=\lambda^2p^r$ and $3D_1s^2-D_2=\pm\lambda^2$}
\end{aligned}\right.\right\}.
\end{align*}
Note that for $\lambda =2$, the condition ``odd" on the prime $p$ should be removed from the above notations. 
\begin{thma}\cite[Theorem 1]{BS01}\label{BST}
Assume that $D_1$ and $D_2$ are coprime positive integers, and $p$ is a prime satisfying $\gcd(D_1 D_2, p)=1$. Then for $\lambda\in \{1, \sqrt{2}, 2\}$, the number of positive integer solutions $(x, y)$ of the Diophantine equation
 \begin{equation}\label{eqd2}
 D_1x^2+D_2=\lambda^2p^y
 \end{equation}
is at most one, except for $$
(\lambda,D_1,D_2,p)\in\Omega:=\left\{\begin{aligned}
&(2,13,3,2),(\sqrt 2,7,11,3),(1,2,1,3),(2,7,1,2),\\
&(\sqrt 2,1,1,5),(\sqrt 2,1,1,13),(2,1,3,7)
\end{aligned}\right\}
$$
and $(D_1, D_2, p)\in
\mathcal{F}\cup \mathcal{G}\cup \mathcal{H}$.
\end{thma}

Note that the authors in \cite{BS01} were unable to determine $(\lambda, D_1, D_2, p)=(2,7,25,2)$  in the set $\Omega$ due to a mild error in calculation. It gives two solutions to \eqref{eqd2}, namely, $(x,y)=(1,3), (17, 9)$. This  comes from simple computation, and it can also be confirmed by \cite{LE93} that these are the only solutions in positive integers corresponding to the above quadruple. 

We also need the following result of the author.
\begin{lem}\cite[Lemma 2.1]{AH20}\label{lemAH20}
 For an integer $k \geq  0$, let $F_k$ (resp. $L_k$) denote the $k$-th Fibonacci (resp. Lucas) number. Then for $\varepsilon = \pm 1$, we have $4F_k - F_{k-2\varepsilon} = L_{k+\varepsilon}$.
\end{lem}
In \cite{CO92}, Cohn completely solved the Diophantine equation $x^2+2^k=y^n$ in positive integers $x, y$ and  $n$, when $k\geq 1$ is an odd integer. We deduce the following lemma from his result (see \cite[Theorem]{CO92}).  
\begin{lem}\label{lem1}
The solutions of the equation 
\begin{equation}\label{eqd3}
x^2+2=3^n, ~~x,y, n\in \mathbb{N}
\end{equation}
are $(x, n)=(1,1), (5,3)$.
\end{lem}
On the other hand, for even positive integer $k$, Arif and Abu Muriefah gave the complete solution of the Diophantine equation $x^2+2^k=y^n$ in positive integers $x, y$ and  $n$. The next lemma can easily be deduced from their result \cite[Theorem 1]{AA01}.  
\begin{lem}\label{lem2}
The solutions of the equation 
\begin{equation}
x^2+4=5^n, ~~x,y, n\in \mathbb{N}
\end{equation}
are $(x, n)=(1,1), (11,3)$.
\end{lem}

\section{Proof of Theorem \ref{thm}}
Assume that $N$ is a generalized Mersenne number such that $N=cx^2$. Then for some prime $p$ and positive integer $n$, we have $M_{p, n}=cx^2$. This can be written as 
\begin{equation}\label{eqp1}
cx^2+p-1=p^n.
\end{equation}
Clearly  $p\nmid cx$ and $\gcd(c, p-1)=1$. It is easy to see from \eqref{eqp1} that there is no generalized Mersenne number of the form $cx^2$ when $2\mid c$.

Let the sets $\Omega, \mathcal{F}, \mathcal{G}$ and $ \mathcal{H}$ be as in Theorem \ref{BST}. If $(c, p-1, p)\in \Omega$, then $(c, p-1, p)=(2,1,3)$ which is not possible.

Now let $(c, p-1, p)\in \mathcal{F}$. Then $\left(F_{k-2\varepsilon}, L_{k+\varepsilon}, F_k\right)=(c, p-1, p)$, and thus by Lemma \ref{lemAH20} $4p-c=p-1$. This implies that $c+p\equiv 1\pmod 4$. 

If $p>2$, then \eqref{eqp1} modulo $4$ gives
$$c\equiv \begin{cases}
1\pmod 4\hspace{7mm} \text{ if } 2\nmid n,\\
2-p\pmod 4 \text{ if } 2\mid n.
\end{cases}$$ 
Thus $c+p\equiv 1\pmod 4$ implies that either $p\equiv 0\pmod 4$ or $2\equiv 1\pmod 4$, and none of these is possible. 
 
 For $p=2$, we have $F_k=2$ and $L_{k+\varepsilon}=1$, which are again not possible.  Therefore,  $(c, p-1, p)\not\in \mathcal{F}$

Assume that $(c, p-1, p)\in \mathcal{G}$. Then $p-1=4p^r-1$ for some positive integer $r$. This implies that $p(4p^{r-1}-1)=0$ and thus $4p^{r-1}=1$  which is not possible. Thus, $(c, p-1, p)\not\in \mathcal{G}$.

Now if $(c, p-1, p)\in \mathcal{H}$, then for odd prime $p$, we have 
\begin{equation}\label{eqp2}
cs^2+p-1=p^r
\end{equation}
and 
\begin{equation}\label{eqp3}
3cs^2-p+1=\pm 1,
\end{equation}
where $r$ and $s$ are positive integers.

From \eqref{eqp3}, we have
$3cs^2=p$ and thus $(c, p,s)=(1,3,1)$. Therefore \eqref{eqp1} becomes
$x^2+2=3^n$, and hence by Lemma \ref{lem1} we have $(x,n)=(1,1), (5, 3)$. This shows that $1$ and $25$ are only generalized Mersenne numbers $M_{3,n}$ which can be written in the form $cx^2$.

Again from \eqref{eqp3}, we have
$3cs^2=p-2$ and thus \eqref{eqp2} gives $p(4-3p^{r-1})=5$, which implies that $(p, r)=(5,1)$. This gives $(c, s)=(1, 1)$ and hence \eqref{eqp1} becomes
$x^2+4=5^n$. By Lemma \ref{lem2}, we get $(x, n)=(1,1), (11, 3)$, which shows that $1$ and $121$ are only generalized Mersenne numbers $M_{5,n}$ that can be written in the form $cx^2$. Thus, we complete the proof by Theorem \ref{BST}. 

\section*{Acknowledgements}
The author would like to thank the anonymous referee for carefully reading this paper.  This work was supported by  SERB-MATRICS (MTR/2021/000762), Govt. of India.


\end{document}